\documentclass{birkjour}
\usepackage{hyperref}
\author{Henryk Trappmann}
\address{Kameruner Str. 9\\13351 Berlin\\ Germany}
\email{henryk@pool.math.tu-berlin.de}
\author{Dimitrii Kouznetsov}
\address{Institute for Laser Science\\University of Electro-Communications
1-5-1 Chofugaoka\\ Chofushi, Tokyo, 182-8585\\ Japan}
\email{dima@uls.uec.ac.jp}

\title
[Uniqueness of Holomorphic Abel Functions]
{Uniqueness of Holomorphic Abel Functions at a Complex Fixed Point Pair}

\begin{document}
\begin{abstract}
  We give a simple uniqueness criterion (and some derived criteria)
  for holomorphic Abel functions and show that Kneser's real analytic
  Abel function of the exponential is subject to this criterion.
\end{abstract}
\subjclass{Primary 30D05} 
\keywords{Abel function; Abel equation; exponential function;
  fractional iterates; holomorphic solution; real analytic}

{This article is a preliminary version of the final article published
in Aequationes Mathematicae with DOI 10.1007/s00010-010-0021-6.
The final publication is available at\\
\url{http://www.springerlink.com/content/u7327835m2850246/}.}

\newpage
\maketitle

\newtheorem{theorem}{Theorem}
\newtheorem{lemma}[theorem]{Lemma}
\newtheorem{proposition}[theorem]{Proposition}
\newtheorem{corollary}[theorem]{Corollary}
\theoremstyle{definition}
\newtheorem{definition}{Definition}
\newtheorem{criterion}{Criterion}

\newcommand \nS {\!\!\!\!\!\!\!\!\!\!\!\!\!\!\!\!\!\!\!}
\newcommand \pS {{~}~{~}}

\newcommand {\iLa}[1]{\label{#1}}         

\newcommand {\iL}[1]{\iLa{#1}} 
\newcommand {\rf}[1]{(\ref{#1})} 

\newcommand \sx {\scalebox}
\newcommand \Lc {L^\ast}
\newcommand{\rme}{\text{e}}
\newcommand \rmi {\text{i}}
\newcommand \sexp {{\rm ksexp}}
\newcommand \slog {{\rm kslog}}
\newcommand \ve {\varepsilon}
\newcommand \bN {\mathbb{N}}
\newcommand \bC {\mathbb{C}}
\newcommand \bR {\mathbb{R}}
\newcommand \bZ {\mathbb{Z}}
\newcommand{\rh}{\varrho}
\newcommand{\wo}{\setminus}
\newcommand{\el}{\ell}

\section{Introduction}
There is a lot of discussion about the ``true'' or ``best'' fractional iterates of
the function $\rme^x$ in the (lay-)mathematical community. In 1949
Kneser \cite{Kneser1949} proved the existence of real analytic fractional
iterates. However Szekeres (a pioneer in developing the theory of fractional iteration
\cite{Szekeres1958}) states 1961 in \cite{Szekeres1961}:

{\it ``The solution of Kneser does not really solve the problem of
`best' fractional iterates of $\rme^x$. Quite apart from practical
difficulties involved in the calculation of Kneser's function on the
real axis, there is no indication whatsoever that the function will
grow more regularly to infinity than any other solution. There is
certainly no uniqueness attached to the solution; in fact if $g(x)$ is
a real analytic function with period 1 and $g'(x)+1>0$ (e.g.
$g(x)=\frac{1}{4\pi}\sin( 2\pi x)$ then 
$B^\ast(x)=B(x)+g\big(B(x)\big)$
 is also
an analytic Abel function of $\rme^x$ which in general yields a different
solution of the equation.''}

A recent discussion with Prof.\ Jean {\'E}calle supports the impression that
no uniqueness criterion was found up today and that there is even evidence
against the existence of a criterion concerned with the
growth-scale or asymptotic behavior at infinity. 

By withdrawing our attention from the purely real analytic behavior
of the Abel function to the behavior in the complex plane we can
succeed in giving a simple uniqueness criterion for Abel functions of
a whole class of real analytic (or 
arbitrary holomorphic) functions with two complex fixed points.

We show the usefulness of the criterion by providing an Abel function that satisfies the criterion.
This is the above mentioned by Kneser constructed Abel
function of $\rme^x$ (which can be easily generalized to functions
$b^x$ with $b>\rme^{1/\rme}$).

We have also a suggestion to numerically compute this Abel function
and the corresponding fractional iterates of $\rme^x$ (also of $b^x$ for
$b>\rme^{1/\rme}$ in generalization) by a method developed in
\cite{Kouznetsov2009}. Several other methods to numerically compute
holomorphic fractional iterates of $\rme^x$ or the holomorphic
Abel function have emerged in the past 
years (for example one is given in \cite{Walker:abelian}). 
A future research goal would be to put them on a thorough
theoretic base (proving convergence and holomorphy) and to verify the
here given uniqueness criterion. 

\section{Motivation}
Our original motivation was the investigation of a fourth stage of
operations after the third stage containing power, exponential and
logarithm.

Different terms for such operations were used in the past like:
``generalized exponential'' and ``generalized logarithm'' by Walker
\cite{Walker:generalized}, ``ultra exponential'' and ``infra
logarithm'' by Hooshmand 
\cite{Hooshmand2006}, ``super-exponential'' by Bromer
\cite{Bromer1987}, tetration and superlogarithm
\cite{Kouznetsov2009}. In this paper we give them the more succinct
names ``4-exponential'' and ``4-logarithm''.

\begin{definition}[4-exponential]
A {\em 4-exponential} to base $b>0$ is a function $f$ that satisfies
\begin{align}
  f(0)&=1\iL{F01}\\
  f(z+1)&=\exp_b\big(f(z)\big)
  \iL{expb}
\end{align}
for all applicable $z$.
\end{definition}
For any $\tilde{f}$ that only satisfies \eqref{expb} and contains 1 in its
codomain: $\tilde{f}(z_0)=1$, the function $f(z)=\tilde{f}(z+z_0)$ is a
4-exponential.
\begin{definition}[4-logarithm]
A {\em 4-logarithm} to base $b>0$ is a
function $g$ that satisfies the Abel equation \eqref{abel} (see
\cite{Kuczma_etal1990}) (for all applicable $z$) with the following
initial condition:
\begin{align}
  g(1)&=0		\iL{g10}\\
  g(\exp_b(z))&=g(z)+1	\iL{abel}.
\end{align}

\end{definition}
For any $\tilde{g}$ that only satisfies \eqref{abel} and has $1$ in
its domain of definition, the function
$g(z)=\tilde{g}(z)-\tilde{g}(1)$ is a 4-logarithm.
Here we set as usual $\exp_b(z)=b^z=\exp\big(\ln(b)z\big)$. 
The inverse of a 4-exponential (if existing) is a 4-logarithm and vice
versa.

On positive integer arguments $z=n$ any 4-exponential $f$ is already determined
to be just the $n$-times application of $\exp_b$ to $1$. 
\begin{align}
f(n)=\exp_{b}^{\circ n}(1)= 
\underbrace{b^{.^{.^{.^b}}}}_{n \times b}\iL{zexpo}
\end{align}
The question however is how to properly extend the function
real and analytic to non-integer arguments. 

The existence of a real analytic strictly
increasing 4-logarithm was proven by Kneser \cite{Kneser1949}. A
non-analytic solution with a uniqueness criterion was given by
Hooshmand in \cite{Hooshmand2006}. A numerical method to compute the real
coefficients of the powerseries development at 0 of a 4-logarithm was
given (though without proof of convergence) by Walker in
\cite{Walker:abelian}. Another numerical method to compute a
real analytic 4-exponential via Cauchy integrals (though also without
convergence proof) was given by Kouznetsov in \cite{Kouznetsov2009}.

A real analytic 4-exponential is expected to have a singularity or
branchpoint at integers $\le -2$ at least on some branch,
because from $f(z+1)=\exp_b(f(z))$ follows
$f(z-1)=\log_b(f(z))$ and by $f(0)=1$ is then $f(-1)=0$
and $f(-2)=\log_b(0)$. To exclude branching we restrict
4-exponentials to
\begin{align}
C_{-2}= \mathbb{C} \backslash  \{x \in \mathbb{R} : x \le -2 \} \iL{C}
\end{align}

It is a conjecture of the authors that holomorphy on the domain
$C_{-2}$ together with $f(z^\ast)=f(z)^\ast$ (complex conjugation) on
$C_{-2}$ implies the uniqueness of the 4-exponential $f$.

From considerations about the uniqueness of
4-logarithms/4-exponentials the following
general uniqueness criterion for Abel functions with two complex
fixed points emerged.

\section{The Uniqueness Criterion}
Before we start we mention some {\em conventions} we use: Usually
curves here are regarded as continuous maps on the open interval
$(-1,1)$. If we however use a curve in a set context then we refer to
the image of the curve, e.g. $\gamma_1\cup \gamma_2 = \gamma_1((-1,1))
\cup \gamma_2((-1,1))$. The disjoint union $C=A\uplus B$ means here
that $C=A\cup B$ and $A\cap B=\emptyset$. The sum $A+z$ of a region
$A\subseteq \bC$ and a number $z\in\bC$ is defined as the region $\{
a+z\colon a\in A\}$. A function being holomorphic 
on a non-open set means that there is a neighborhood of each point of
the set where the function is holomorphic. $\log_b(z)=\log(z)/\log(b)$
means the principal branch $-\pi<\Im(\log(z))\le\pi$ of the
logarithm if not stated otherwise. ``Continuable'', ``continuation''
and ``continue'' always refer to {\em analytic} continuation.

\begin{definition}[Abel function, initial
  curve/region]
We call a function $\alpha$, holomorphic on $D$, an {\em Abel function
  of $F$} iff it satisfies the Abel equation
\begin{align}
  \alpha(F(z))&=\alpha(z)+1
  \iL{eq:abel}
\end{align}
for all $z\in D\cap F^{-1}(D)$. $F$ is sometimes called the
{\em base function}.

A curve $\gamma\colon (-1,1)\to\bC$ on which $F$ is
holomorphic is called an {\em initial curve of $F$} iff $\gamma$ and
$F\circ\gamma$ are injective
and disjoint and $\gamma(-1)\neq \gamma(1)$ are two fixed points of
$F$. (To be explicit define $\gamma(\pm 1):=\lim_{t\to\pm 1} \gamma(t)$.)

Under these conditions $\gamma\cup (F\circ \gamma) \cup
\{\gamma(-1),\gamma(1)\}$ is a closed Jordan
curve. We call its inner (bounded) component $C$ joined with $\gamma$
and $F\circ\gamma$ the {\em initial region} of $\gamma$ denoted by
$I_F(\gamma):=\gamma\uplus C\uplus F(\gamma)=\overline{C}\setminus
\{\gamma(-1),\gamma(1)\}$. 
%
\iL{def:initial}
\end{definition}

As a side remark: the following theorem has a slight similarity to the criterion given
in \cite{contreras:2007:remarks}. Both criteria include
injectivity. However the criterion given in 
\cite{contreras:2007:remarks} applies to a function with one parabolic
fixed point (moreover a self-map of the unit disk of hyperbolic step zero),
while the criterion here applies to a function with two fixed points
(which can not occur for a holomorphic self-map of some simply
connected region).

\begin{theorem}
Let $\gamma$ be an initial curve of $F$, let $H$ be its initial region and $d\in H$.
There is at most one function $\alpha$ satisfying Criterion
\ref{criterion:A}. 
\iL{theorem}
\end{theorem}
\begin{criterion}
  The function $\alpha$ is an on $H$ holomorphic and injective Abel function of
  $F$, $\alpha(d)=0$ and $\bigcup_{k\in\bZ} \left(\alpha(H)+k\right)=\bC$.
\iL{criterion:A}
\end{criterion}
\begin{proof} 
Assume there are two such Abel functions $\alpha_1\colon H \leftrightarrow
T_1$ and $\alpha_2\colon H\leftrightarrow T_2$  holomorphic and injective on
$H$. For the rest of this proof we write $\alpha_j$ when referring to $\alpha_1$ as
well as to $\alpha_2$. The inverse function $\alpha_j^{-1}\colon T_j\leftrightarrow H$
satisfies:
\begin{align*}
  \alpha_j^{-1}(z+1)&=F\left(\alpha_j^{-1}(z)\right)
\end{align*}
for all $z$ such that $z,z+1\in T_j$.
So we have two biholomorphic functions 
\begin{align*}
  q_1&:=\alpha_2\circ \alpha_1^{-1}\colon T_1\leftrightarrow T_2 &
  q_2:=\alpha_1\circ \alpha_2^{-1}\colon T_2\leftrightarrow T_1
\end{align*}
with the property
\begin{align*}
   q_1(z+1)
  &=\alpha_2\left(\alpha_1^{-1}(z+1)\right)=\alpha_2\left(F\left(\alpha_1^{-1}(z)\right)\right)\\
  &=\alpha_2\left(\alpha_1^{-1}(z)\right)+1 = q_1(z)+1
\end{align*}
for each $z$ with $z,z+1\in T_1$; and generally
\begin{align}
   q_j(z+1)= q_j(z)+1
\iL{eq:delta}
\end{align}
for each $z$ with $z,z+1\in T_j$.

We define $q_{j,k}\colon T_j+k\to\mathbb{C}$ by $q_{j,k}(z+k)=q_j(z)+k$.  By our
property $q_j(z+1)=q_j(z)+1$ the function
$q_{j,k}$ and $q_{j,k+1}$ coincide on the intersection
$(T_j+k)\cap (T_j+k+1)$ which contains the curve $\alpha\circ\gamma + k+1$. In
conclusion $q_j$ can be continued to
the whole complex plane. So lets consider $q_j$ to
be an entire function.

For $z\in T_1$ we have $ q_1(z) = q_2^{-1}(z)$. But the only entire
functions that have an entire inverse are linear functions. By the
values $ q_1(0)=0$ and $ q_1(1)=1$ it can only be the
identity. So $\alpha_2(\alpha_1^{-1}(z))=z$ for $z\in T_1$ and hence $\alpha_2=\alpha_1$
on $H$.
\end{proof}

Now one may argue that the Abel function may depend on the initial
region. This is not the case as shown by the next theorem.

%
\begin{theorem}
  Let $D$ be a connected set with $d\in D$ and let $F$ be holomorphic
  on $D$, there can be at most one on $D$ holomorphic Abel function 
  $\alpha$ of $F$ with the property that $\alpha'(d)\neq 0$ and that $\alpha$
  satisfies Criterion \ref{criterion:A} on some initial region
  $H\subset D$.
\end{theorem}
\begin{proof}
  Assume that there were two such Abel functions $\alpha_1$ and
  $\alpha_2$ such that each $\alpha_j$ is
  holomorphic on $D$ and satisfies Criterion \ref{criterion:A} on the
  initial region $H_j$.  
  We follow the proof of Theorem \ref{theorem}, listing only the
  modifications.

  We consider biholomorphic $\alpha_j\colon H_j\leftrightarrow T_j$ and
  holomorphic $q_j\colon T_j\to\bC$ and achieve \eqref{eq:delta}.
  In conclusion each $q_j$ is entire.
  The function $z\mapsto q_1(\alpha_1(z))-\alpha_2(z)$ is holomorphic
  on $D$. It is constantly 0 on $H_1$ and hence on 
  $D$. That's why $q_1(0)=q_1(\alpha_1(d))=\alpha_2(d)=0$. Vice versa
  $q_2(0)=0$. 

  There is a neighborhood $V$ of $d$ where $\alpha_1$ and $\alpha_2$ is
  injective. Hence $\alpha_j^{-1}$, $q_1=\alpha_2\circ \alpha_1^{-1}$
  and $q_2=\alpha_1\circ \alpha_2^{-1}$ are injective on
  $V'=\alpha_1(V)\cap \alpha_2(V)\ni 0$. But then $q_1=q_2^{-1}$ on
  $V'\cap q_2(V')$.
\end{proof}

Now we want to make the criterion $\bigcup_{k\in\bZ}
\left(\alpha(H)+k\right)=\bC$ a bit more accessible and show that it is a
consequence of $\lim_{t\to\pm 1}\Im(\alpha(\gamma(t)))=\pm\infty$.
Before we start with the actual proof, we need a little
insight into how curves divide the complex plane. We know by the
Jordan curve theorem that each simple closed curve divides the sphere
into two simply connected components. Considering the sphere $\bC\cup
\{\infty\}$ we know that each injective curve $\zeta\colon (-1,1)\to
\bC$
with $\lim_{t\to\pm 1} \zeta(t) = \infty$ divides the complex plane
into two parts; where $\infty$ is the complex infinity and $\lim_{t\to
  \pm 1} \zeta(t) =\infty$ means that for each $r>0$ there is a $t_1$
and $t_0$ such
that $\left|\zeta(t)\right| > r$ for all $t>t_1$ and all $t<t_0$.

A particular subclass of such curves are the injective
curves $\zeta$ with $\lim_{t\to \pm 1}\Im(\zeta(t)) = \pm\infty$. Here
$\infty$ is the
real infinity and $\lim_{t\to \pm 1}\Im(\zeta(t)) = \pm\infty$ means that for
each $u>0$ there are $t_0,\; t_1\in (-1,1)$ such that $\Im(\zeta(t)) >
u$ for all $t>t_1$ and $\Im(\zeta(t))<-u$ for all $t<t_0$.

\begin{definition}[left/right component/ray]
For an injective curve $\zeta\colon (-1,1)\to\bC$ we call a component {\em
  left} (resp.\ {\em right}) if  
it contains a left (resp.\ right) ray, where a left (resp.\ right) ray
is a set of the form $\{ z_0\mp x \colon x> 0\}$ for some $z_0\in \bC$.
\end{definition}
The following lemma shows that these properties indeed behave as
expected.
\begin{lemma}[$L$, $R$]
  Let $\zeta\colon (-1,1)\to \bC$ be an injective curve with $\lim_{t\to\pm
    1} \Im(\zeta(t)) = \pm \infty$ and let $P$ and $Q$ be the two
  components the plane is divided into; $\bC =  P\uplus \zeta
  \uplus Q$. 

  1) Then either $P$ is left and $Q$ is right or vice versa. We denote the left
  (resp.\ right) component with $L(\zeta)$ (resp.\ $R(\zeta)$).
  Moreover there exist left (resp.\ right) rays contained in the left
  (resp.\ right) component for each prescribed imaginary part $y$.

  2) If $\zeta+d$ is disjoint from $\zeta$ then $\overline{R}(\zeta+d)\subset
  R(\zeta)$ in the case $d>0$ and $\overline{L}(\zeta+d)\subset L(\zeta)$ in the
  case $d<0$. Here the overline means the closure of the component
  (it is equal to the union of the component with $\zeta$).
\iL{lemma}
\end{lemma}
\begin{proof}
  For part 1) of the lemma we show that for each prescribed
  imaginary part $y$ there exists a left ray contained in one component
  and a right ray contained in the other component. We show further
  that the union of all left (resp.\ right) rays such that each is contained in
  either $P$ or $Q$ must be contained in either $P$ or $Q$, which
  ensures the ``either'' in part 1).

  We consider the indices of the intersection of the horizontal line
  $Y=\{x+iy\colon x\in\bR\}$ with the curve $\zeta$,
  $T=\{t\colon \zeta(t)\in Y\}$. Now let $t_0=\inf(T)$ and $t_1=\sup(
  T)$; neither can $t_0=-1$ nor $t_1=1$ because in this case there
  would be a sequence of $t\to\pm 1$ such that $\Im(\zeta(t)) = y$ in
  contradiction to $\Im(\zeta(t))\to \pm \infty$.  

  Hence $-1<t_0\le t_1<1$ and $\Im(\zeta(t))>y$ for all
  $t>t_1$ and $\Im(\zeta(t))<y$ for all $t<t_0$. Let
  $x_0=\Re(\zeta(t_0))$ and $x_1=\Re(\zeta(t_1))$ then it is clear that
  $Y_0:=\{x+iy\colon x<x_0\}$ is completely contained in a component as
  well as $Y_1:=\{x+iy\colon x>x_1\}$ is completely contained in a
  component. The compound $Y_0\cup
  \zeta([x_0,x_1])\cup Y_1$ divides the plane into an upper half $H_1$
  which is divided by $\zeta((x_1,1))$ and a lower part $H_0$ which
  is divided by $\zeta((-1,x_0))$. 

  If there was a path $\beta\colon [0,1]\to\bC$ that connects
  $Y_0$ with $Y_1$, i.e.\ $\beta(0)\in Y_0$ and $\beta(1)\in Y_1$ then
  we choose $s_0=\beta^{-1}(\sup(\Re(\beta\cap Y_0))+\rmi y)$ and
  $s_1=\beta^{-1}(\inf(\Re(\beta\cap Y_1))+\rmi y)$.
  The restriction of $\beta$ to the
  non-empty interval $(s_0,s_1)$ still connects $Y_0$ to $Y_1$ but
  does neither intersect $Y_0$ nor $Y_1$. As the path is also not 
  allowed to intersect $\zeta([t_0,t_1])$ it must either be contained
  in $H_0$ or in $H_1$ but then it would intersect $\zeta((-1,x_0))$
  or $\zeta((x_1,1))$. 

  If we have two left rays with the imaginary parts $y_1$ and $y_2$ then
  we can do the above construction of $t_0(y)$ for 
  $y=y_1,y_2$. Without restriction let $t_0(y_1)<t_0(y_2)$, then 
  $\zeta([t_0(y_1),t_0(y_2)])$ has a minimum of the real part and we
  can connect both left rays by a vertical line with a smaller real
  part. 

  Now to part 2: We consider $d>0$. First it is easy to see
  that one point of $\zeta+d$ lies in $R(\zeta)$ and hence the whole
  $\zeta+d$ is contained in $R(\zeta)$. Vice versa $\zeta$ must be
  contained in $L(\zeta+d)$. So $\zeta$ does not intersect
  $R(\zeta+d)$ and we conclude that $R(\zeta+d)$ must be either contained in
  $L(\zeta)$ or in $R(\zeta)$. But $R(\zeta)$ has common points with
  $R(\zeta+d)$ as a right ray $Y$ is contained in $R(\zeta)$ and the
  ray $Y+d$ is contained in $R(\zeta+d)$.
\end{proof}

\begin{criterion}
  Under the preconditions of Theorem \ref{theorem},
  and $\gamma$ being rectifiable, the following criterion implies
  Criterion \ref{criterion:A}: 

  The curve $\zeta = \alpha\circ\gamma$ is injective, $\zeta \cap
  (\zeta+1)=\emptyset$ and $\lim_{t\to\pm
    1}\Im(\zeta(t))=\pm\infty$; where $\alpha$ is an on $H$ holomorphic Abel function of $F$ with $\alpha(d)=0$.
\iL{criterion:B}
\end{criterion}
\begin{proof}
The first two conditions of the criterion state that $\alpha$ is
injective on $\gamma\cup (F\circ \gamma)$ because $\alpha\circ F\circ
\gamma=\zeta+1$ by \eqref{eq:abel}.
By a theorem about univalent functions (see \cite{Markushevich1965}
Theorem 4.8) the injectiveness of $\alpha$ on the rectifiable boundary with
finitely many exceptions (e.g.\ $\gamma(-1)$ and $\gamma(1)$) implies
the injectiveness on the enclosed region (inclusive its boundary)
which is $H$.

Now we show that $\bigcup_{k\in\bZ} \left(\alpha(H)+k\right)=\bC$:
The image $\alpha(H)$ is bounded by $\zeta$ and $\zeta+1$ and it must
be simply connected, that's why 
\begin{align*}
  \alpha(H) &= \overline{R}(\zeta) \cap \overline{L}(\zeta+1)
\end{align*}
If we now unite consecutive pieces $\alpha(H)+k$, $k\in\bZ$, we get with
$\overline{R}(\zeta-1)\supset \overline{R}(\zeta)$ and
$\overline{L}(\zeta)\subset \overline{L}(\zeta+1)$ from Lemma \ref{lemma} that:
\begin{align*}
  \bigcup_{k=k_0}^{k_1} \left(\alpha(H)+k\right) &=
  \overline{R}(\zeta+k_0)\cap \overline{L}(\zeta+1+k_1)
\end{align*}

By part 1) of Lemma \ref{lemma} we also know that for each
$z=x+i y$ there is a (negative) $k_0\in\bZ$ such that $z\in R(\zeta+k_0)$ (right
ray with starting point $x-k_0+i y$ is contained in $R(\zeta)$)
and a (positive) $k_1\in\bZ$ such that $z\in L(\zeta+k_1)$ (left ray with
starting point $x-k_1+i y$ is contained in $L(\zeta)$). So 
$\alpha(H)$ translated by all integers cover the
whole complex plane:
\begin{align}
  \bigcup_{k\in\bZ} \left(\alpha(H)+k\right) &= \bC 
\iL{eq:C}
\end{align}
\end{proof}

\begin{criterion}
  Under the preconditions of Theorem \ref{theorem}
  the following criterion implies Criterion \ref{criterion:B}: 

  The real function $f(t)=\Im(\alpha(\gamma(t)))$ is strictly
  increasing and $\lim_{t\to\pm 1}f(t) = \pm \infty$; where $\alpha$
  is an on $H$ holomorphic Abel function of $F$ with $\alpha(d)=0$.
\iL{criterion:C}
\end{criterion}
\begin{proof}
  First it is clear that $\zeta:=\alpha\circ\gamma$ is injective as no imaginary
  value can be taken twice.
  Further the correlation $\mu$ given by
  $\Im(\zeta(t))\mapsto \Re(\zeta(t))$, $t\in (-1,1)$, is a welldefined
  function on $\bR$.
  Then $\mu+1>\mu$ is a function not intersecting $\mu$, and hence $\zeta+1$
  does not intersect $\zeta$.
\end{proof}


\section{Application to Kneser's construction}
In this section we apply the uniqueness Criterion \ref{criterion:A} to the 
4-logarithm (i.e.\ Abel function of $\exp_b$) $\Psi$ as constructed by
Kneser \cite{Kneser1949}.

\begin{figure}
\begin{center}
\sx{.32}{\begin{picture}(600,460)
\put(20,0){\includegraphics{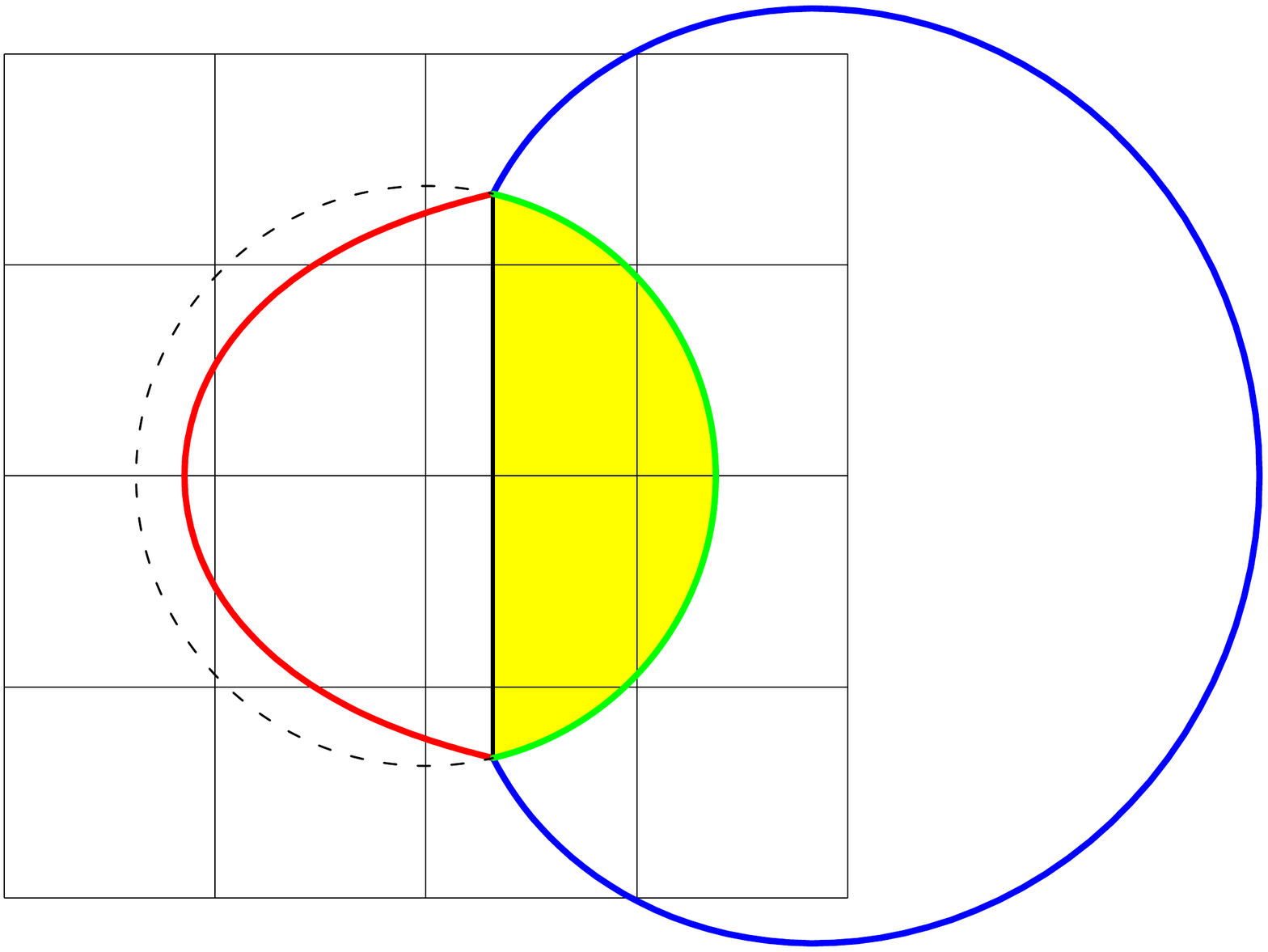}}
\put(204,452){\sx{2.5}{$\Im(z)$}}
\put(204,422){\sx{2.5}{$2$}}
\put(204,322){\sx{2.5}{$1$}}
\put(204,222){\sx{2.5}{$0$}}
\put(182,122){\sx{2.5}{$-1$}}
\put(182, 22){\sx{2.5}{$-2$}}
\put(  2,208){\sx{2.5}{$-2$}}
\put(102,208){\sx{2.5}{$-1$}}
\put(317,208){\sx{2.5}{$1$}}
\put(417,208){\sx{2.5}{$2$}}
\put(457,208){\sx{2.5}{$\Re(z)$}}
\put(140,260){\sx{2.5}{$\log_b(H)$}}
\put(280,250){\sx{4}{$H$}}
\put(400,300){\sx{4}{$\exp_b(H)$}}
\put(110,170){\sx{2.5}{$\log_b(\ell)$}}
\put(243,160){\sx{5}{$\ell$}}
\put(320,140){\sx{5}{$b^\ell$}}
\put(520, 40){\sx{6}{$b^{b^\ell}$}}
\put(233,352){\sx{4.4}{$L$}}
\put(233, 82){\sx{4.4}{$\Lc$}}
\end{picture}}
\end{center}
\caption{Contours $\log_b(\ell)$, $\ell$ and $b^\ell$ in the complex $z$-plane,
for base $b=\rme$. \iL{fighen}}
\end{figure}

Related to Kneser's construction we use the choice 
\begin{align}
  H=\{z\in \bC\colon \Re(z)\ge \Re(L),\; |z|\le |L|\}\wo \{L,\Lc\}
\iL{H}
\end{align}
as initial region which is depicted in Figure \ref{fighen}
for $b=\rme$ (which is the only base that Kneser considered).
Here $L$ is the fixed point of $\log_b$ in the upper half plane.
The straight line 
\begin{align}
  \ell(t)=\Re(L)+\rmi\Im(L)t, \quad -1<t<1,
\iL{ell}
\end{align} between $L$ and its
complex conjugate $\Lc$ is the left boundary of $H$ and $b^{\ell}$ is
the right boundary of $H$. 

\begin{lemma}
  If $b>e^{1/e}$ then $\ell$ in \eqref{ell} is an initial curve and
  hence $H$ is an initial region of $\exp_b$.
\end{lemma}
\begin{proof}
  We show that $b^\ell$ is injective and does not intersect $\ell$. By $b^L=L$ we
  know that $b^{\Re(L)}=|b^L|=|L|$ and hence 
  \begin{align*}
    b^{\el(t)}=b^{\Re(L)+\rmi\Im(L)t}=|L|\rme^{\rmi\Im(L)\ln(b)t}
  \end{align*}
  which is an arc with radius $|L|$ centered in 0 (shown with a dashed line)
  starting at angle $-\Im(L)\ln(b)$ and ending at angle $\Im(L)\ln(b)$. 
  This is true for any non-real conjugated fixed
  point pair of $\exp_b$. For $b>\rme^{1/\rme}$ there are no real fixed points
  of $\exp_b$ and the fixed point pair of $\log_b$, which is the one closest to the real axis, has
  $\pm\Im(L) \in (-\pi/\ln(b),\pi/\ln(b))$. This assures that
  $b^{\ell}$ does not overlap itself, i.e.\ that it is injective.
%
\end{proof}
Let us --- without proof --- enumerate some counterexamples of initial
curves:
$\log_b\circ \ell$ is initial for $e^{1/e}<b<e^{\pi/2}$ but $\ell$ has
zero or negative real part for $b\ge e^{\pi/2}$ and is hence no more contained in
the default domain of the logarithm. The curve $b^\ell$ is
initial for $b=\rme$ but there are bigger bases where $b^{b^\ell}$ intersects
itself. Each of $\rme^\ell$, $\rme^{\rme^\ell}$,
$\rme^{\rme^{\rme^\ell}}$ are injective but
$\rme^{\rme^{\rme^{\rme^\ell}}}$ is not, so the first two are initial
for $\exp$. The straight line connecting any other conjugated fixed point
pair is not initial because the image under $\exp_b$ is a circle with
radius $|L|$ (winding at least once around 0).

\newcommand{\Ho}{\mathfrak{H}_0}
To get familiar with Kneser's construction we recapitulate the main
steps he does in \cite{Kneser1949} with a slight generalization to bases
$b>\rme^{1/\rme}$. He starts with the K{\oe}nigs function $\chi$ of
$\exp_b$ at the fixed point $L$ and shows that it can be continued to
nearly the whole upper half plane $\mathfrak{H}=\{z\in \bC\colon \Re(z)\ge
0\}\wo\{0,1,b,b^b,b^{b^b},\dots\}$. He also shows that it is injective
on $\mathfrak{H}$ because its inverse can be continued to an entire
function.

Particularly $\chi$ is injective on the ``half'' initial region
\begin{align*}
\Ho&=\{z\in \bC\colon \Re(z)\ge \Re(L),\; \Im(z)\ge 0,\;|z|\le |L|\}\wo
\{L,1\}.
\end{align*} 
Note that Kneser initially works with $\Ho$
containing $L$ while he later switches to consider $\Ho$ without $L$
which is the definition we use above. $\Ho$ is related to our initial region $H$ via $\Ho \cup \Ho^\ast  = H \wo \{1\}$. 

The K{\oe}nigs function satisfies the Schr\"oder equation
$\chi(b^z)=c\chi(z)$ on $z,b^z\in \mathfrak{H}$ where 
$c=\exp_b'(L)=\log(L)$. Next he sets $\psi(z)=\log(\chi(z))$ for a
suitable region of the logarithm which has the cut outside $\chi(
\Ho
)$. It
satisfies $\psi(b^z)=\psi(z)+\log(c)$ on $z\in 
\Ho$. Last he conformally maps the union of
$\psi(\mathfrak{H}_0) + k\cdot\log(c)$, $k\in\bZ$, via
the Riemann mapping theorem to the upper halfplane, say with the
conformal map $\rh$. 

The resulting function $\Psi = \rh \circ\log\circ \chi$ is at least
defined on $
\Ho
$ and satisfies
$\Psi(b^z)=\Psi(z)+1$ on $
\Ho$. It can even be continued to $z=
1
$. It is real analytic on $\bR\cap 
\Ho 
$ and can hence analytically continued to the conjugate
region $
\Ho
^\ast$ (and also to the whole real line).

From his construction the following things are
important for us. First:
$\Psi$ is injective and holomorphic on $\Ho\cup\{1\}$ and the
image is contained in the upper halfplane. Hence the conjugate
continuation to $H$ is also injective and holomorphic. Second: 
the by integer translated regions $\Psi(\Ho\cup\{1\})$ cover
the whole upper halfplane and hence the by integer translated regions
$\Psi(H)$ cover the whole complex plane. 

The application of Theorem \ref{theorem} gives:
\begin{theorem}
  The by Kneser in \cite{Kneser1949} constructed real analytic 4-logarithm
  $\alpha=\Psi$ and each generalization to base $b>\rme^{1/\rme}$
  satisfies Criterion \ref{criterion:A} on $H$ given in \eqref{H} 
  (where $d=1$ and $F(z)=b^z$).
\end{theorem}

 
\section{About Computation of the
  4-Logarithm/4-Exponential}
As Kneser uses the Riemann mapping theorem in his construction,
it is very difficult to approach computationally. Instead we use
here a rather fast numerical method given in \cite{Kouznetsov2009}
utilizing the Cauchy integral formula to compute a 4-exponential which
we denote with $\sexp_b$ and its inverse 
with $\slog_b$. This method is originally described for
$b=\rme$ but can be extended to arbitrary bases
$b>\rme^{1/\rme}$. Unfortunately it lacks a proof of convergence (but
if it converges then $\slog_b$ satisfies the Abel equation and is
holomorphic). It
is the conjecture of the authors that this method does converge and
that it satisfies the uniqueness Criterion \ref{criterion:C}. 

\begin{figure}
\begin{center}
\sx{3}{\begin{picture}(50,66)
\put(0,0){\includegraphics{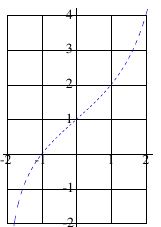}} 
\put(0,0){\includegraphics{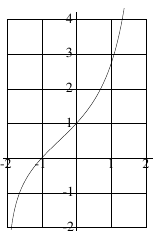}} 
\put(10,64.5){\sx{.34}{$\sexp_{b}(x)$}}
\put(33,65){\sx{.33}{$b=\rme$}}
\put(42,63){\sx{.33}{$b=2$}}
\put(45,20){\sx{.4}{$x$}}
\end{picture}}
\end{center}
\caption{
4-exponential $\sexp_b$ to base  $b=\rme$ (thick solid) and $b=2$
(dashed)
on the real axis.} 
\iL{fig01}
\begin{center}
\sx{1.8}{\begin{picture}(140,140)
\put(8,8){\includegraphics{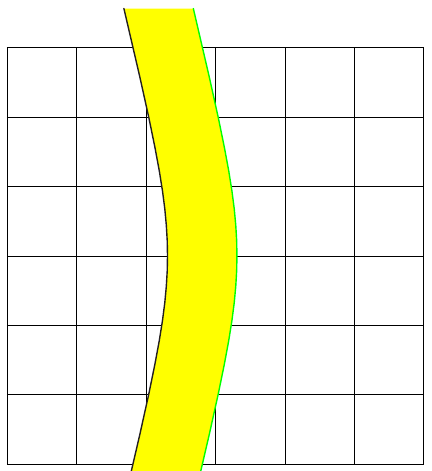}}
\put(8,8){\includegraphics{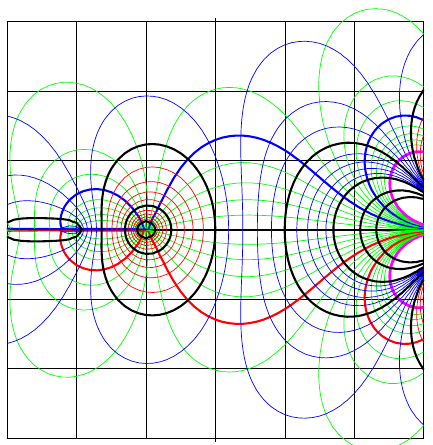}}
\put(50,116){\sx{.44}{$\slog(H)$}}
\put(-1,133){\sx{.7}{$\Im(z)$}}
\put(4,108){\sx{.6}{$2$}}
\put(4, 88){\sx{.6}{$1$}}
\put(4, 68){\sx{.6}{$0$}}
\put( 0, 48){\sx{.6}{$-1$}}
\put( 0, 28){\sx{.6}{$-2$}}
\put(23, 2){\sx{.6}{$-2$}}
\put(43, 2){\sx{.6}{$-1$}}
\put(69, 2){\sx{.6}{$0$}}
\put(89, 2){\sx{.6}{$1$}}
\put(109,2){\sx{.6}{$2$}}
\put(120,1){\sx{.7}{$\Re(z)$}}
\end{picture}}
\end{center}
\caption{Contour plot of the function $\sexp_\rme(z)$ showing
  lines of
  constant modulus and of constant phase. The shaded (yellow) region is
  $\slog_\rme(H)$. 
Its right boundary is $\slog_\rme\left(\rme^{\ell}\right)$ which
corresponds to $|\sexp_\rme(z)|=|L|$. \iL{figmapa3} 
}
\end{figure}
For real values of the argument, this 4-exponential 
is plotted in Figure \ref{fig01} for $b=\rme,2$.
The plot of $\slog_\rme$ in the complex plane (depicted in figure
\ref{figmapa3}) hints towards Criterion \ref{criterion:C}. 

\subsection{The Fractional Iterates of the Exponential}
The combination of $\sexp$ and $\slog$ allows to define
fractional/continuum iterative powers of the exponential via (see
e.g.\ \cite{Szekeres1958} or \cite{Kuczma_etal1990}):
\begin{align}
\exp_{b}^{\circ c}(z)=\sexp_{b}\big(c+\slog_{b}(z)\big)
\iL{expbcz}
\end{align}
Unlike regular iterates at a hyperbolic fixed point which are analytic at
the fixed point, the above iterates of $\exp_b$ have
a branch point at both complex fixed points. Care must be taken to
determine a principal branch/region/cut of 
$\slog_b$ which is then also the principal region of the iterative
power. It should be chosen such that $c+\slog_{b}(z)\in C_{-2}$.
The map of the function $\exp^{\circ 0.5}$ in the complex plane is plotted in reference \cite{Kouznetsov2010}, 
showing a behavior similar to the iterative square root of the factorial
$\sqrt{~!~}$ computed in \cite{Kouznetsov2010}. 

\begin{figure}
\begin{center}
\sx{1.2}
{\begin{picture}(190,190)
\put(6,6){\includegraphics{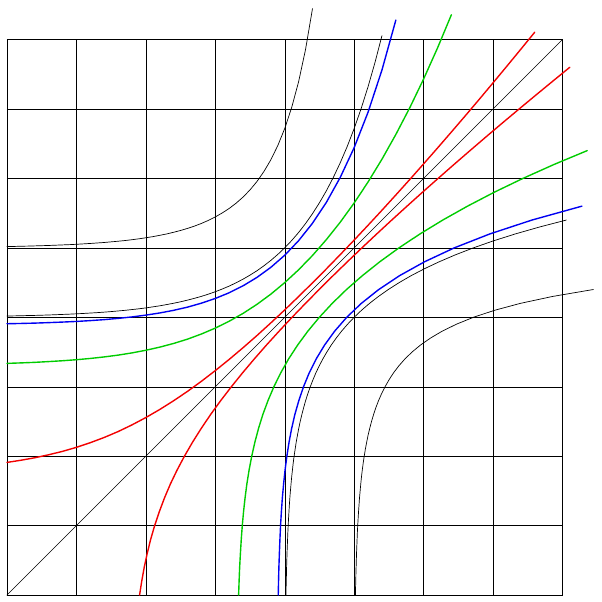}}
\put(86,177){\sx{.66}{$c=2$}}
\put(102,169.4){\sx{.66}{$c=1$}}
\put(108,176){\sx{.66}{$c=0.9$}}
\put(134,176){\sx{.66}{$c=0.5$}}
\put(143,169.4){\sx{.66}{$c=0.1$}}
\put(169,170){\sx{.66}{$c=0$}}
\put(170,160){\sx{.66}{$c=-0.1$}}
\put(170,137){\sx{.66}{$c=-0.5$}}
\put(170,121){\sx{.66}{$c=-0.9$}}
\put(170,114){\sx{.66}{$c=-1$}}
\put(171,97){\sx{.66}{$c=-2$}}
\put(-1,171){\sx{.8}{$y=\exp^{\circ c}(x)$}}
\put(1,145){\sx{.8}{3}}
\put(1,125){\sx{.8}{2}}
\put(1,105){\sx{.8}{1}}
\put(1, 85){\sx{.8}{0}}
\put(0, 65){\sx{.8}{-1}}
\put(0, 45){\sx{.8}{-2}}
\put(0, 25){\sx{.8}{-3}}
\put(24, 0){\sx{.8}{-3}}
\put(44, 0){\sx{.8}{-2}}
\put(64, 0){\sx{.8}{-1}}
\put(86, 0){\sx{.8}{0}}
\put(106,0){\sx{.8}{1}}
\put(126,0){\sx{.8}{2}}
\put(146,0){\sx{.8}{3}}
\put(170,0){\sx{.8}{$x$}}
\end{picture}}
\end{center}
\caption{ Function 
$y=\exp^{\circ c}(x)$ calculated by equation (\ref{expbcz}) for
$c=0$, $\pm 0.1$, $\pm 0.5$, $\pm 0.9$, $\pm 1$, $\pm 2$ versus $x$. }
\iL{fig14}
\end{figure}

For several {\em real} $c$ we plot $\exp^{\circ c}(x)$ versus
$x$ in Figure \ref{fig14}. For $c=1$ it is indeed the usual exponential,
for $c=0$ it is the identity function, and for $c=-1$ it is the
logarithm.

\bibliographystyle{amsplain}
\bibliography{main}

\providecommand{\bysame}{\leavevmode\hbox to3em{\hrulefill}\thinspace}
\providecommand{\MR}{\relax\ifhmode\unskip\space\fi MR }
\providecommand{\MRhref}[2]{%
  \href{http://www.ams.org/mathscinet-getitem?mr=#1}{#2}
}
\providecommand{\href}[2]{#2}
\begin{thebibliography}{10}

\bibitem{Bromer1987}
N.~Bromer, \emph{Superexponentiation}, Math. Mag. \textbf{60} (1987), no.~3,
  169--174.

\bibitem{contreras:2007:remarks}
M.D. Contreras, S.D. Madrigal, and C.~Pommerenke, \emph{Some remarks on the
  {A}bel equation in the unit disk}, J. Lond. Math. Soc., II. Ser. \textbf{75}
  (2007), no.~3, 623--634.

\bibitem{Hooshmand2006}
M.H. Hooshmand, \emph{Ultra power and ultra exponential functions}, Integral
  Transforms Spec. Funct. \textbf{17} (2006), no.~8, 549--558.

\bibitem{Kneser1949}
H.~Kneser, \emph{Reelle analytische {L}{\"o}sungen der {G}leichung
  {$\varphi(\varphi(x))=e^x$} und verwandter {F}unktionalgleichungen}, J. Reine
  Angew. Math. \textbf{187} (1949), 56--67.

\bibitem{Kouznetsov2009}
D.~Kouznetsov, \emph{Solution of $f(x+1)=\exp(f(x))$ in complex $z$-plane},
  Math. Comp. \textbf{78} (2009), 1647--1670.

\bibitem{Kouznetsov2010}
D.~Kouznetsov and H.~Trappmann, \emph{Superfunctions and square root of
  factorial}, Mosc. Univ. Phys. Bull. \textbf{65} (2010), no.~1, 6--12.

\bibitem{Kuczma_etal1990}
M.~Kuczma, B.~Choczewski, and R.~Ger, \emph{Iterative functional equations},
  Cambridge University Press, 1990.

\bibitem{Markushevich1965}
A.~I. Markushevich, \emph{Theory of functions of a complex variable},
  Prentice-Hall, 1965.

\bibitem{Szekeres1958}
G.~Szekeres, \emph{Regular iteration of real and complex functions}, Acta Math.
  \textbf{100} (1958), 203--258.

\bibitem{Szekeres1961}
\bysame, \emph{Fractional iteration of exponentially growing functions}, J.
  Austral. Math. Soc. \textbf{2} (1961), 301--320.

\bibitem{Walker:generalized}
Peter~L. Walker, \emph{Infinitely differentiable generalized logarithmic and
  exponential functions}, Math. Comput. \textbf{57} (1991), no.~196, 723--733.

\bibitem{Walker:abelian}
\bysame, \emph{On the solutions of an {A}belian functional equation}, J. Math.
  Anal. Appl. \textbf{155} (1991), no.~1, 93--110.

\end{thebibliography}

\end{document}